\documentclass[11pt,twoside,letterpaper]{article}
\usepackage{pslatex}
\usepackage{fancyhdr}
\usepackage[english,francais]{babel}
\usepackage{graphicx}
\usepackage{geometry}

\usepackage{amsmath}
\usepackage{amssymb}
\usepackage{amsfonts}
\usepackage{amsthm,amscd}

\newtheorem{theoreme}{Th\'eor\`eme}[section]
\newtheorem{lemme}[theoreme]{Lemme}
\newtheorem{proposition}[theoreme]{Proposition}
\newtheorem{corollaire}[theoreme]{Corollaire}
\newtheorem{definition}[theoreme]{D\'efinition}
\newtheorem{remarque}[theoreme]{Remarque}

\newtheorem{probleme}[theoreme]{Probl\`eme}

\theoremstyle{remark}

\numberwithin{equation}{section}

\begin{document}
\vskip 0.4in
\title{\bfseries\scshape{Les tournois (-1)-critiques}}
\author{Houmem Belkhechine\thanks{E-mail address: houmem@gmail.com}\\Facult\'e des Sciences de Gab\`es, Cit\'e Riadh, Zirig 6072
Gab\`es, Tunisie \and Imed Boudabbous\thanks{E-mail adress:
imed.boudabbous@gmail.com}\\Institut Pr\'eparatoire aux \'Etudes
d'Ing\'enieurs de Sfax, route Menzel Chaker Km 0.5 - 3018
Sfax,Tunisie \and Jamel Dammak\thanks{E-mail:
jdammak@yahoo.fr}\\Facult\'e des Sciences de Sfax, BP 802, 3018
Sfax, Tunisie}

\date{}
\maketitle \thispagestyle{empty} \setcounter{page}{1}

% ------- [First Page Running Head] - place it immediately after title! ------
\thispagestyle{fancy} \fancyhead{} \fancyhead[L]{{\bf
Communications in Mathematical Analysis}\\
Volume 2, Number 2, (2007)\\ ISSN \ 0000-00} % put \label{lastpage-xx} on the last page!
\fancyhead[R]{\sl www.commun-math-anal.org\\} \fancyfoot{}
\renewcommand{\headrulewidth}{0pt}

%\noindent{\bf Abstract}

\begin{abstract}
\selectlanguage{francais}
 \noindent \'Etant donn\'e
un tournoi $T= ( S,A)$, une partie $X$ de $S$ est un intervalle de
$T$ lorsque pour tous $a, b\in X$ et $ x\in S-X$, $(a,x)\in A$ si
et seulement si $(b,x)\in A$. Par exemple, $\emptyset$,
$\{x\}(x\in S)$ et $S$ sont des intervalles de $T$, appel\'es
intervalles triviaux. Un tournoi, dont tous les intervalles sont
triviaux, est ind\'ecomposable; sinon, il est d\'ecomposable. Un
sommet $x$ d'un tournoi ind\'ecomposable $T$ est critique si le
tournoi $T- x$ est d\'ecomposable. En 1993, J.H. Schmerl et W.T.
Trotter ont caract\'eris\'e les tournois dont tous les sommets
sont critiques, appel\'es tournois critiques. Ces tournois ont un
cardinal impair $\geq 5$. Pour chaque entier impair $m \geq 5$, il
existe trois tournois critiques de cardinal $m$. Dans cet article,
nous caract\'erisons les tournois qui admettent un unique sommet
non critique, que nous appelons tournois (-1)-critiques. Ces
tournois ont un cardinal impair $\geq 7$. Pour chaque entier
impair $m \geq 7$, il existe $3m - 15$ tournois $(-1)$-critiques
de cardinal $m$.

\end{abstract}

%\noindent {\bf AMS Subject Classification:} 62G05; 62G20.

\vspace{.08in} \noindent \textbf{Mots cl\'es}: Critique, Graphe
d'ind\'ecomposabilit\'e, Intervalle, Tournoi Ind\'ecomposable.

\selectlanguage{english}
\begin{abstract}
\noindent {\bf The (-1)-critical tournaments. } Given a tournament
$T=(V,A)$, a subset $X$ of $V$ is an interval of $T$ provided that
for any $a, b\in X$ and $ x\in V-X$, $(a,x)\in A$ if and only if
$(b,x)\in A$. For example, $\emptyset$, $\{x\}(x\in V)$ and $V$
are intervals of $T$, called trivial intervals. A tournament, all
the intervals of which are trivial, is indecomposable; otherwise,
it is decomposable. A vertex $x$ of an indecomposable tournament
is critical if $T-x$ is decomposable. In 1993, J.H. Schmerl and
W.T. Trotter characterized the tournaments, all the vertices of
which are critical, called critical tournaments. The cardinality
of these tournaments is odd. Given an odd integer $m \geq 5$,
there exist three critical tournaments of cardinality $m$. and
there are exactly three critical tournaments for each such a
cardinality. In this article, we characterize the tournaments
which admit a single non critical vertex, that we call
(-1)-critical tournaments. The cardinality of these tournaments is
odd. Given an odd integer $m \geq 7$, there exist $3m-15$
(-1)-critical tournaments of cardinality $m$.
\end{abstract}
%\selectlanguage{french}
\section{Introduction} %

Un {\it graphe} ({\it orient\'e}) $G = (S(G), A(G))$ ou $(S, G)$,
est constitu\'e d'un ensemble fini $S$ de sommets et d'un ensemble
$A$ de couples de sommets distincts, appel\'es {\it arcs} de $G$.
L'\emph{ordre} (ou le {\it cardinal}) du graphe $G$ est le nombre
de ses sommets. \`A chaque partie $X$ de $S$ est associ\'e le {\it
sous-graphe} $G(X) = (X, A \cap (X \times X))$ de $G$ induit par
$X$. Pour $X \subseteq S$ (resp. $x \in S$), le graphe $G(S-X)$,
o\`u $S-X = \{s \in S : s \notin X\}$, (resp. $G(S-\{x\}$) est
not\'e $G-X$ (resp. $G-x$). \'Etant donn\'es deux graphes $G = (S,
A)$ et $G' = (S, A')$, une bijection $f$ de $S$ sur $S'$ est un
{\it isomorphisme} de $G$ sur $G'$ si pour tous $x$, $y \in S$,
$(x, y) \in A$ si et seulement si $(f(x), f(y)) \in A'$. Lorsqu'un
tel isomorphisme existe, on dit que $G$ et $G'$ sont {\it
isomorphes}, et on note $G \simeq G'$.

Un graphe $G = (S, A)$ est un {\it tournoi} lorsque pour tous $x
\neq y \in S$, on a: $(x, y) \in A$ si et seulement si $(y, x)
\notin A$. \'Etant donn\'e un tournoi \ $T=(S,A)$, pour tous
sommets distincts $x$, $y$ de $S$, la notation $x \longrightarrow
y$ signifie $(x, y) \in A$, et on dit, dans ce cas, que $x$ domine
$y$. Pour toutes parties disjointes $I$ et $J$ de $S$, on note $I
\longrightarrow J$ lorsque pour tout $(x, y) \in I \times J$, $x
\longrightarrow y$. La notation $I \sim J$ signifie que $I
\longrightarrow J$ ou $J \longrightarrow I$. De m\^{e}me, pour
tout $x \in S$ et pour tout $Y \subseteq S-\{x\}$, $x
\longrightarrow Y$ (resp. $Y \longrightarrow x$) signifie $x
\longrightarrow y$ (resp. $y \longrightarrow x$) pour tout $y \in
Y$. La notation $x \sim Y$ signifie que $x \longrightarrow Y$ ou
$Y \longrightarrow x$. Pour tout $x \in S$, on pose $V_{T}^{-}(x)=
\{y \in S$ : $y \longrightarrow x \}$ et $V_{T}^{+}(x)= \{y \in S$
: $x \longrightarrow y \}$. On introduit une relation
d'\'equivalence, not\'ee $\equiv$, sur les couples de sommets
distincts de $T$, d\'efinie comme suit: $(x, y) \equiv (u, v)$ si
$(x, y) = (u, v)$ ou $\mid \{(x, y), (u, v)\} \cap A \mid \neq 1$.

Un tournoi $T$ est un {\it ordre total} (ou une {\it cha\^ine}, ou
une {\it liste}), lorsque pour tous $x, y, z\in S(T) $, si $x
\longrightarrow y$ et $y \longrightarrow z$, alors $x
\longrightarrow z$. Un ordre total d'ordre $k$ est aussi appel\'e
\emph{$k$-cha\^ine}. Pour deux sommets distincts $a$ et $b$ d'un
ordre total $T$, $a < b$ signifie $a \longrightarrow b$. La
notation $T = a_{0} < \cdots < a_{n}$ signifie que $T$ est l'ordre
total d\'efini sur $S=\{a_{0}, \cdots, a_{n}\}$ par $A(T) =
\{(a_{i},a_{j}) : i < j\}$. L'ordre total usuel $0 < \cdots < n$
est not\'e $L_{n+1}$.

 \`A tout tournoi $T=(S, A)$ est associ\'e son
tournoi {\it dual} $T^{\star}=(S, A^{\star})$, o\`u $A^{\star} =
\{(x, y) : (y,x) \in A\}$.

\'Etant donn\'e un tournoi $T=(S, A)$, une partie $I$ de $S$ est
un \emph{intervalle} [4, 5, 7] (ou un {\it clan} \cite{E}) de $T$
lorsque pour $a, b\in I$ et $ x\in S-I$, $(a, x) \equiv (b, x)$.
Par exemple, $\emptyset $, $\{x\}$ o\`u $x \in S$, et $S$ sont des
intervalles de $T$, appel\'es les intervalles {\it triviaux} de
$T$. Un tournoi est {\it ind\'ecomposable} [5, 7] (ou {\it
primitif} \cite{E}) si tous ses intervalles sont triviaux et il
est {\it d\'ecomposable} dans le cas contraire.

Un sommet $x$ d'un tournoi ind\'ecomposable $T$ est dit {\it
critique} si le tournoi $T-x$ est d\'ecomposable. Soit $T$ un
tournoi ind\'ecomposable \`a au moins $5$ sommets. Le tournoi $T$
est dit {\it critique} si tous ses sommets sont critiques. On
g\'en\'eralise cette d\'efinition en disant que le tournoi $T$ est
{\it $(-k)$-critique} lorsqu'il admet exactement $k$ sommets non
critiques. Afin de rappeler la caract\'erisation des tournois
critiques, nous introduisons, pour tout entier $n \geq 1$, les
tournois $T_{2n+1}$, $U_{2n+1}$ et $V_{2n+1} $ d\'efinis sur $\{0,
\cdots, 2n\}$ comme suit.

\begin{enumerate}

\item $A(T_{2n+1}) = \{(i,j) : j-i \in \{1, \cdots, n\}\
mod. \ 2n+1\}$.
\item $U_{2n+1}$ est le tournoi obtenu \`a partir de l'ordre total
$L_{2n+1}$ en inversant les arcs reliant deux sommets pairs, de
sorte que $A(U_{2n+1}) = \{(i,j) : i<j$ et $i$ ou $j$ est impair
$\} \cup \{(i,j) : i>j $ et $i$ et $j$ sont pairs$\}$.
\item $V_{2n+1}(\{0,\cdots, 2n-1\}) = 0<\cdots<2n-1$ et $V_{V_{2n+1}}^{+}(2n) = \{2i : i \in \{0,
\cdots,n-1\}\}$.

\end{enumerate}
Remarquons que $T_{3}=U_{3}=V_{3}=\{\{0, 1, 2\}, \{(0, 1), (1, 2),
(2, 0) \}\}$.
\begin{proposition} {{\rm (J. H. Schmerl et W. T. Trotter \cite{ST})}}
\`A un isomorphisme pr\`es, les tournois critiques sont les
tournois $T_{2n+1}$,  $U_{2n+1}$ et $V_{2n+1}$, o\`u $n \geq 2$.
\end{proposition}

Dans cet article, nous caract\'erisons les tournois
$(-1)$-critiques, r\'epondant ainsi, dans le cas des tournois, \`a
une question pos\'ee par Y. Boudabbous et P. Ille \cite{BI}.
Contrairement aux tournois $T_{2n+1}$, les tournois $U_{2n+1}$ et
$V_{2n+1}$ apparaissent dans la morphologie de ces tournois que
nous pr\'esentons \`a partir de leur unique sommet non critique. A
cet effet, nous d\'efinissons pour tout entier $n \geq 3$ et pour
tout entier $k \in \{1, \cdots, n-2\}$, les tournois $
E_{2n+1}^{2k+1} $, $F_{2n+1}^{2k+1}$, $G_{2n+1}^{2k+1}$ et
$H_{2n+1}^{2k+1}$ d\'efinis sur $\{0, \cdots, 2n\}$ comme suit.
\begin{enumerate}

\item $E_{2n+1}^{2k+1}(V_{E_{2n+1}^{2k+1}}^{-}(2k+1)) = L_{2k+1}$, $E_{2n+1}^{2k+1}(V_{E_{2n+1}^{2k+1}}^{+}(2k+1)) = 2k+2 < \cdots <
2n$ et pour tout $(x, y) \in V_{E_{2n+1}^{2k+1}}^{+}(2k+1) \times
V_{E_{2n+1}^{2k+1}}^{-}(2k+1)$, $x \longrightarrow y$ si et
seulement si $x$ et $y$ sont pairs.

\item $F_{2n+1}^{2k+1}(V_{F_{2n+1}^{2k+1}}^{-}(2k+1)) = U_{2k+1}$, $F_{2n+1}^{2k+1}(V_{F_{2n+1}^{2k+1}}^{+}(2k+1)) = 2k+2 < \cdots <
2n$ et pour tout $(x, y) \in V_{F_{2n+1}^{2k+1}}^{+}(2k+1) \times
V_{F_{2n+1}^{2k+1}}^{-}(2k+1)$, $x \longrightarrow y$ si et
seulement si $x$ et $y$ sont pairs.

\item $G_{2n+1}^{2k+1}(V_{G_{2n+1}^{2k+1}}^{-}(2k+1)) = U_{2k+1}$, $G_{2n+1}^{2k+1}(V_{G_{2n+1}^{2k+1}}^{+}(2k+1)) \simeq V_{2n-2k-1}$
avec $2k+2 < \cdots < 2n-1$ et pour tout $(x, y) \in
V_{G_{2n+1}^{2k+1}}^{+}(2k+1) \times
V_{G_{2n+1}^{2k+1}}^{-}(2k+1)$, $x \longrightarrow y$ si et
seulement si $x =2n$ et $y$ est pair.

\item $H_{2n+1}^{2k+1}(V_{H_{2n+1}^{2k+1}}^{-}(2k+1)) = V_{2k+1}$, $H_{2n+1}^{2k+1}(V_{H_{2n+1}^{2k+1}}^{+}(2k+1)) \simeq
V_{2n-2k-1}$ avec $2k+2 < \cdots < 2n-1$ et pour tout $(x, y) \in
V_{H_{2n+1}^{2k+1}}^{+}(2k+1) \times
V_{H_{2n+1}^{2k+1}}^{-}(2k+1)$, $x \longrightarrow y$ si et
seulement si $x= 2n$ et $y = 2k$.

\end{enumerate}

Observons que dans les tournois $T = E_{2n+1}^{2k+1}$,
$F_{2n+1}^{2k+1}$, $G_{2n+1}^{2k+1}$ ou $H_{2n+1}^{2k+1}$
d\'efinis ci-dessus, chacun des sous-tournois $T(V_{T}^{-}(2k+1))$
et $T(V_{T}^{+}(2k+1))$ est ou bien une cha\^ine, ou bien un
tournoi critique non isomorphe \`a un tournoi de la classe
$\{T_{2p+1}\}_{p \geq 2}$. Les cardinaux respectifs de ces
sous-tournois sont $2k+1$ et $2n-2k-1$.

Pour tout entier $n \geq 3$, on d\'esigne par $\mathcal{E}_{2n+1}$
(resp. $\mathcal{F}_{2n+1}$, $\mathcal{F}_{2n+1}^{\star}$,
$\mathcal{G}_{2n+1}$, $\mathcal{G}_{2n+1}^{\star}$,
$\mathcal{H}_{2n+1}$), la classe des $n-2$ tournois
$\{E_{2n+1}^{2k+1}\}_{1 \leq k \leq n-2}$ (resp.
$\{F_{2n+1}^{2k+1}\}_{1 \leq k \leq n-2}$,
$\{(F_{2n+1}^{2k+1})^{\star}\}_{1 \leq k \leq n-2}$ ,
$\{G_{2n+1}^{2k+1}\}_{1 \leq k \leq n-2}$,
$\{(G_{2n+1}^{2k+1})^{\star}\}_{1 \leq k \leq n-2}$,
$\{H_{2n+1}^{2k+1}\}_{1 \leq k \leq n-2}$).

Remarquons alors le fait suivant.
\begin{remarque}

\'Etant donn\'e un entier $n \geq 3$, si $T$ est un tournoi de la
classe $\mathcal{E}_{2n+1}$ (resp. $\mathcal{H}_{2n+1}$), alors
$T^{\star}$ est aussi un tournoi de la classe $\mathcal{E}_{2n+1}$
(resp. $\mathcal{H}_{2n+1}$)).
\end{remarque}
\noindent{\em Preuve\/. } Il suffit de remarquer que pour $k \in
\{1, \cdots, n-2\}$, la permutation $\sigma$ de $\{0, \cdots,
2n\}$ d\'efinie par : pour tout $q \in \{0, \cdots, 2n \}$,
$\sigma(q) = 2n-q$ (resp. $\sigma(q) = 2n-q-1$ si $q \in \{0,
\cdots, 2n\} - \{2n, 2k, 2k+1\}$, $\sigma(2n) = 2(n-k-1)$,
$\sigma(2k) = 2n$ et $\sigma(2k+1) = 2(n-k-1) +1$), est un
isomorphisme de $(E_{2n+1}^{2k+1})^{\star}$ (resp.
$(H_{2n+1}^{2k+1})^{\star}$) sur $E_{2n+1}^{2(n-k-1) +1}$ (resp.
$H_{2n+1}^{2(n-k-1) +1}$).{\hspace*{\fill}$\Box$\medskip}

La caract\'erisation suivante des tournois $(-1)$-critiques, est
le principal r\'esultat de cet article.

\begin{theoreme}\`A un isomorphisme pr\`es, les tournois $(-1)$-critiques sont les
tournois $E_{2n+1}^{2k+1}$, $F_{2n+1}^{2k+1}$,
$(F_{2n+1}^{2k+1})^{\star}$, $G_{2n+1}^{2k+1}$,
$(G_{2n+1}^{2k+1})^{\star}$ et $H_{2n+1}^{2k+1}$, o\`u $n \geq 3$
et $1 \leq k \leq n-2$. De plus, le sommet $2k+1$ est l'unique
sommet non critique de chacun de ces tournois.

\end{theoreme}
%\clearpage
%\begin{figure}[h]
%\begin{center}

%\includegraphics[width=14cm]{figs.jpg}
%\end{center}
%\end{figure}
%Les figures ci-dessus repr\'esentent, \`a un isomorphisme pr\`es
%et au dual pr\`es, les quatre tournois $(-1)$-critiques \`a $7$
%sommets. Les arcs manquants sont orient\'es du haut vers le bas.
%Le sommet 3 est l'unique sommet non critique de chacun de ces
%tournois.

\section{Tournois ind\'ecomposables}

\begin{definition}

Soit $T = (S, A)$ un tournoi. \`A toute partie $X$ de $S$ telle
que $\mid X \mid \geq 3$ et le sous-tournoi $T(X)$ est
ind\'ecomposable, on associe les parties de $S-X$ suivantes.
\begin{itemize}
\item $[X] = \{x \in S-X: x \sim X\}$.

\item Pour tout $u \in X$, $X(u) = \{x \in S-X: \ \{u, x\}$ est un intervalle de
$T(X \cup \{x\})\}$.

\item $Ext(X) = \{x \in S-X: \ T(X \cup \{x\})$ est ind\'ecomposable$\}$.
\end{itemize}

\end{definition}

Rappelons le lemme suivant.

\begin{lemme} {{\rm(A. Ehrenfeucht et G. Rozenberg \cite{E})}}
\label{Lem.(Ehrenfeucht+Rozenberg)} %Lemme
Soient $T = (S, A)$ un tournoi et $X$ une partie de $S$ tels que
$\mid X \mid \geq 3$ et $T(X)$ est ind\'ecomposable. La famille
$\{X(u): \ u \in X\} \cup \{Ext(X), [X]\}$ forme une partition de
$S-X$ . De plus, les assertions suivantes sont v\'erifi\'ees.

\begin{itemize}

\item Soient $u \in X$, $x \in X(u)$ et $y \in S-(X \cup X(u))$. Si
$T(X \cup \{x, y\})$ est d\'ecomposable, alors $\{u, x\}$ est un
intervalle de $T(X \cup \{x, y\})$.

\item Soient $x \in [X]$ et  $y \in S-(X \cup [X])$. Si
$T(X \cup \{x, y\})$ est d\'ecomposable, alors $X \cup \{y\}$ est
un intervalle de $T(X \cup \{x, y\})$.

\item Soient $x \neq y \in Ext(X)$. Si
$T(X \cup \{x, y\})$ est d\'ecomposable, alors $\{x, y\}$ est un
intervalle de $T(X \cup \{x, y\})$.

\end{itemize}

\end{lemme}

De ce lemme d\'ecoule le r\'esultat suivant.

\begin{corollaire} {{\rm(A. Ehrenfeucht et G. Rozenberg \cite{E})}}

Soit $T = (S, A)$ un tournoi ind\'ecomposable. Si $X$ est une
partie  de $S$ telle que $\mid X \mid \geq 3$, $\mid S-X \mid \geq
2$ et $T(X)$ est ind\'ecomposable, alors il existe deux sommets
distincts $x$ et $ y$ de $S-X$ tels que $T(X\cup\{x, y\})$ est
ind\'ecomposable.

\end{corollaire}

\section{Graphe d'ind\'ecomposabilit\'e}

Rappelons d'abord qu'un graphe {\it symetrique} (ou {\it non
orient\'e}) est un graphe $G$ tel que pour tous $x \neq y \in
S(G)$, on a : $(x, y) \in A(G)$ si et seulement si $(y, x) \in
A(G)$. On consid\`ere alors que, dans un tel graphe $G$, $A(G)$
est un ensemble de paires de sommets distincts de $S(G)$,
appel\'ees {\it ar\^etes} de $G$. Par exemple, le {\it chemin}
$P_{n}$ de longueur $n-1$ et le {\it cycle} $C_{n}$ de longueur $n
\geq 3$ sont les graphes non orient\'es d\'efinis sur $\{0,
\cdots, n-1\}$ de la fa\c{c}on suivante. Pour tous $i$, $j \in
\{0, \cdots, n-1\}$, $\{i, j\}$ est une ar\^ete de $P_{n}$ si
$\mid i- j\mid= 1$. Le cycle $C_{n}$ est alors obtenu \`a partir
de $P_{n}$ en ajoutant l'ar\^ete $\{0, n-1\}$. Tout graphe
isomorphe \`a $P_{n}$ (resp. $C_{n}$) est appel\'e chemin (resp.
cycle). Une relation d'\'equivalence $\mathcal{R}$ est d\'efinie
sur $S(G)$ comme suit. Pour tous $x \neq y \in S(G)$, $x$
$\mathcal{R}$ $y$ s'il existe une suite $x_{0} = x, \cdots,
x_{n}=y $ de sommets de $G$ telle que pour tout $i \in \{0,
\cdots, n-1\}$, $\{x_{i}, x_{i+1}\} \in A(G)$. les classes
d'\'equivalence de $\mathcal{R}$ sont appel\'ees {\it composantes
connexes} de $G$. Pour tout sommet $x \in S(G)$, on pose
$V_{G}(x)= \{y \in S$: $\{x, y \} \in A(G)\}$. Lorsque $V_{G}(x) =
\emptyset$, on dit que $x$ est un sommet {\it isol\'e} de $G$.

La notion de {\it graphe d'ind\'ecomposablit\'e} \`a \'et\'e
introduite par P. Ille [1, 5] de la fa\c{c}on suivante. \`A chaque
tournoi $T = (S, A)$ est associ\'e son graphe
d'ind\'ecomposabilit\'e $I(T)$ d\'efini sur $S$ comme suit. Pour
tous $x \neq y \in S$, $\{x, y\}$ est une ar\^ete de $I(T)$ si $T
- \{x, y\}$ est ind\'ecomposable. Ce graphe est un outil important
dans notre construction des tournois $(-1)$-critiques.

Notons qu'un tournoi $T$ et son dual $T^{\star}$ ont les m\^emes
intervalles. Il s'ensuit que $T$ et $T^{\star}$ ont les m\^emes
sommets critiques ainsi que le m\^eme graphe
d'ind\'ecomposabilit\'e.

Dans la suite de ce paragraphe, nous \'etudions le graphe
d'ind\'ecomposabilit\'e d'un tournoi $(-1)$-critique.

Rappelons, d'adord, les deux lemmes suivants.

\begin{lemme} {{\rm(Y. Boudabbous et P. Ille \cite{BI})}} Soient
$T = (S, A)$ un tournoi ind\'ecomposable et $x$ un sommet critique
de $T$. Alors $\mid V_{I(T)}(x) \mid \leq 2$ et on a:

\begin{itemize}
\item Si $V_{I(T)}(x) = \{y\}$, o\`u $y \in S$, alors $T-\{x, y\}$
est un intervalle de $T-x$.
\item Si $V_{I(T)}(x) = \{y, z\}$, o\`u $y \neq z \in S$, alors $\{y,
z\}$ est un intervalle de $T-x$.

\end{itemize}

\end{lemme}

\begin{lemme}{{\rm(Y. Boudabbous et P. Ille \cite{BI})}}
Le graphe d'ind\'ecomposabilit\'e d'un tournoi $(-1)$-critique
admet une unique composante connexe de cardinal $\geq 2$.

\end{lemme}

Le lemme suivant pr\'ecise l'ordre d'un tournoi $(-1)$-critique
\begin{lemme}
L'ordre d'un tournoi $(-1)$-critique est impair et sup\'erieur ou
\'egal \`a 7.
\end{lemme}
\noindent{\em Preuve\/. }

Les tournois \`a 4 sommets sont, \`a un isomorphisme pr\`es, au
nombre de quatre et sont tous d\'ecomposables. Il s'ensuit que les
tournois ind\'ecomposables \`a 5 sommets sont critiques. Ainsi, il
n'existe aucun tournoi $(-1)$-critique d'ordre 5.

 Soit $T$ un tournoi ind\'ecomposable \`a au moins 3 sommets.
Pour tout sommet $x$ de $T$, il existe deux sommets $y \neq z$ de
$T-x$ tels que $T(\{x, y, z\}) \simeq U_{3}$. En effet, autrement,
il existe un sommet $\alpha$ de $T$ tel que $V_{T}^{-}(\alpha)
\longrightarrow V_{T}^{+}(\alpha)$. Si $\mid V_{T}^{-}(\alpha)
\mid = \mid V_{T}^{+}(\alpha) \mid =1$ alors $T \simeq L_{3}$, ce
qui contredit l'ind\'ecomposabilit\'e de $T$. Sinon,
$V_{T}^{-}(\alpha)$ ou $V_{T}^{+}(\alpha)$ est un intervalle non
trivial de $T$, une contradiction. Supposons, \`a pr\'esent, que
le tournoi $T$ est $(-1)$-critique et d\'esignons par $a$ son
unique sommet non critique. D'apr\`es ce qui pr\'ec\`ede, il
existe deux sommets $b \neq c$ de $T-a$ tel que $T(\{a, b, c\})
\simeq U_{3}$. Si le tournoi $T$ est d'ordre pair, alors par une
suite finie d'applications du corollaire 2.3, on obtient un sommet
$\omega \in S(T)-\{a\}$ tel que le tournoi $T-\omega$ est
ind\'ecomposable. Contradiction.

{\hspace*{\fill}$\Box$\medskip}

Le r\'esultat suivant compl\`ete le lemme 3.1 dans le cas des
tournois $(-1)$-critiques.

\begin{lemme} Le sommet non critique $a$ d'un tournoi
$(-1)$-critique $T = (S, A)$ est tel que $\mid V_{I(T)}(a) \mid =
2$.
\end{lemme}

\noindent{\em Preuve\/. } Comme $T-a$ est un tournoi
ind\'ecomposable qui est, d'apr\`es le lemme 3.3, d'ordre pair, il
s'ensuit que $T-a$ n'est ni critique ni $(-1)$-critique. Il existe
alors deux sommets distincts $x$, $y \in S-\{a\}$ tels que les
tournois $T-\{a, x\}$ et $T-\{a, y\}$ sont ind\'ecomposables.
Ainsi $\{x, y\} \subseteq V_{I(T)}(a)$, de sorte que $\mid
V_{I(T)}(a) \mid \geq 2$. Supposons que $\mid V_{I(T)}(a) \mid
\geq 3$ et consid\'erons 3 sommets deux \`a deux distincts $x$,
$y$ et $z$ de $V_{I(T)}(a)$. On pose $X = S - \{a, x\}$, $Y = S -
\{a, y\}$ et $Z = S - \{a, z\}$. Les tournois $T$, $T(X)$, $T(Y)$
et $T(Z)$ sont ind\'ecomposables avec, d'apr\`es le lemme 3.3,
$\mid X \mid = \mid Y \mid = \mid Z \mid \geq 5$. Les tournois
$T-x$, $T-y$ et $T-z$ \'etant d\'ecomposables, alors $a \notin
Ext(X) \cup Ext(Y) \cup Ext(Z)$. Nous montrons d'abord qu'il
existe $(u, v, w) \in X \times Y \times Z$ tel que $a \in X(u)
\cap Y(v) \cap Z(w)$. Supposons, par exemple, qu'il n'existe pas
un $u \in X$ tel que $a \in X(u)$, ce qui \'equivaut \`a dire,
d'apr\`es le lemme 2.2, que $a \in [X]$. Quitte \`a remplacer $T$
par $T^{\star}$, on peut supposer que $a \longrightarrow X$ et
donc $x \longrightarrow a$. Remarquons que $a \notin [Y]$, sinon,
comme $z \in X \cap Y$ et $a \longrightarrow X$, alors $a
\longrightarrow Y$ et en particulier $a \longrightarrow x$, une
contradiction. De m\^eme $a \notin [Z]$. Comme de plus, $a \notin
Ext(Y) \cup Ext(Z)$, alors, d'apr\`es le lemme 2.2, il existe $(v,
w) \in Y \times Z$ tel que $a \in Y(v) \cap Z(w)$.

On a forc\'ement $v \neq w$. Sinon, comme $a \longrightarrow X$ et
$a \in Y(v)$ (resp. $a \in Z(v)$), alors $v \longrightarrow S -
\{v, a, x, y\}$ (resp. $v \longrightarrow S - \{v, a, x, z\}$). Il
s'ensuit que $v \longrightarrow S - \{v, a, x\}$ et en particulier
$v \neq x$, autrement $\{x, a\} \longrightarrow X$, ce qui
contredit l'ind\'ecomposabilit\'e de $T$. Comme de plus $\mid X
\mid \geq 5$, alors $S - \{v, a, x\}$ est un intervalle non
trivial du tournoi ind\'ecomposable $T(X)$. Contradiction.

De plus, $\{v, w\} \cap \{x, y, z\} \neq \emptyset$. Sinon, comme
d'une part $a \longrightarrow X$, en particulier $a
\longrightarrow \{v, w\}$, et d'autre part $a \in Y(v)$ (resp. $a
\in Z(w) $), il s'ensuit que $v \longrightarrow w$ (resp. $w
\longrightarrow v$). Contradiction. Cela nous am\`ene \`a
distinguer les cas suivants.

\begin{itemize}
\item $x \notin \{v, w\}$. Dans ce cas $\{v, w\} \cap \{y, z\} \neq
\emptyset$. Si, par exemple, $v=z$, alors $\{a, z\}$ et $\{a, w\}$
sont des intervalles respectifs de $T-y$ et de $T-z$, de sorte que
$\{z, w\}$ est un intervalle non trivial du tournoi
ind\'ecomposable $T-\{a, y\}$. Contradiction.

\item $x \in \{v, w\}$. Supposons, par exemple, que $v=x$. Dans ce cas $a \in
Y(x)$, et comme $x \longrightarrow a \longrightarrow X$, alors $x
\longrightarrow S-\{x, y\}$. Il s'ensuit que $S- \{a, x, y\}$ est
un intervalle non trivial du tournoi ind\'ecomposable $T - \{a,
y\}$. Contradiction.

\end{itemize}

Il est, \`a pr\'esent, \'etabli qu'il existe $(u, v, w) \in X
\times Y \times Z$ tel que $a \in X(u) \cap Y(v) \cap Z(w)$. D'une
part, les sommets $u$, $v$ et $w$ sont deux \`a deux distincts;
autrement, si par exemple $u = v$, alors $\{a, u\}$ est un
intervalle de chacun des tournois $T-x$ et $T-y$ ce qui implique
que $\{a, u\}$ est un intervalle non trivial du tournoi
ind\'ecomposable $T$. Contradiction. D'autre part, $\{u, v, w\}
\cap \{x, y, z\} \neq \emptyset$; sinon, pour tout $\alpha \in
\{u, v, w \}$, $\{u, v, w\} - \{\alpha\}$ est un intervalle de
$T(\{u, v, w \})$. Une contradiction, car dans chacun des deux
tournois \`a 3 sommets ($U_{3}$ et $L_{3}$), il existe une paire
de sommets qui n'est pas un intervalle. Supposons alors, par
exemple, que $u = y$. Les paires $\{u, a\}$ et $\{v, a\}$ sont
alors des intervalles respectifs des tournois $T-x$ et $T-u$.
Ainsi $\{u, v \}$ est un intervalle non trivial, de $T - \{a, x\}$
si $v \neq x$, de $T-a$ si $v = x$. Une contradiction, puisque les
tournois $T-\{a, x\}$ et $T-a$ sont ind\'ecomposables.
{\hspace*{\fill}$\Box$\medskip}

Le lemme ci-dessus, permet de pr\'eciser la composante connexe,
non r\'eduite \`a un singleton, du graphe d'ind\'ecomposabilit\'e
d'un tournoi $(-1)$-critique.

\begin{lemme} Si $T$ est un tournoi (-1)-critique, alors $I(T)(\mathcal{C})$
est un chemin, o\`u $\mathcal{C}$ est l'unique composante connexe
de $I(T)$ qui n'est pas r\'eduite \`a un singleton.

\end{lemme}

\noindent{\em Preuve\/. } Rappelons que, d'apr\`es le lemme 3.2,
$I(T)$ admet une unique composante connexe $\mathcal{C}$, non
r\'eduite \`a un singleton. De plus, d'apr\`es les lemmes 3.1 et
3.4, pour tout sommet $x \in \mathcal{C}$, on a: $1 \leq \mid
V_{I(T)(\mathcal{C})}(x) \mid \leq 2$. Il s'ensuit que
$I(T)(\mathcal{C})$ est ou bien un chemin ou bien un cycle.
Supposons que $I(T)(\mathcal{C})$ est le cycle $C_{n}$ d\'efini
sur $\mathbb{Z}/n\mathbb{Z} = \{0, \cdots, n-1\}$ o\`u $n \geq 3$
et o\`u 0 est l'unique sommet non critique de $T$. Soit $i \in
\{1, \cdots, n-1 \}$. Comme $i$ est un sommet critique de $T$ avec
$V_{I(T)}(i) = \{i-1, i+1\}$, alors $\{i-1, i+1\}$ est un
intervalle de $T-i$, sans \^etre un intervalle de $T$, de sorte
que $(i-1, i) \equiv (i, i+1)$. Quitte \`a remplacer $T$ par
$T^{\star}$, on peut supposer que $0 \longrightarrow 1$. Il
s'ensuit que pour tout $i \in \{0, \cdots, n-1\}$, $i
\longrightarrow i+1$, en particulier $n-1 \longrightarrow 0$. Si
$j > 1$ est un sommet impair de $\mathcal{C}$, alors $j-1$ est un
sommet critique de $T$ avec $V_{I(T)}(j-1) = \{j-2, j\}$. Il
s'ensuit que $\{j-2, j\}$ est un intervalle du tournoi
$T-\{j-1\}$, et comme $j-1 \neq 0$, alors $(0, j) \equiv (0,
j-2)$. Ainsi, pour tout sommet impair $k$ de $\mathcal{C}$, on a:
$(0, 1) \equiv (0, 3) \equiv \cdots \equiv (0, k)$, et comme $0
\longrightarrow 1$, alors $0 \longrightarrow k$. Il s'ensuit que
$0 \longrightarrow \{k \in \mathcal{C} : k$ est impair\}. Comme
$n-1 \longrightarrow 0$, alors $n$ est impair. Posons $n-1 = 2m$
et distinguons les deux cas suivants.
\begin{itemize}
\item $\mathcal{C} \neq S$. Soit $x \in S - \mathcal{C}$. Comme
pour tout $i \in \{1, \cdots, 2m\}$, $\{i-1, i+1\}$ est un
intervalle de $T-i$, alors $(x, 0) \equiv (x, 2) \equiv \cdots
\equiv (x, 2m)$ et $(x, 0) \equiv (x, 2m-1) \equiv \cdots \equiv
(x, 1)$. Il s'ensuit que $x \sim \mathcal{C}$ et donc
$\mathcal{C}$ est un intervalle non trivial du tournoi
ind\'ecomposable $T$. Contradiction.
\item $\mathcal{C} = S$. D'une part, $0 \longrightarrow 1$ donc $1 \longrightarrow \{i
\in S-\{0\} : i$ est pair\}; d'autre part, $2m-1 \longrightarrow
2m$ donc $\{i \in S : i$ est impair\} $\longrightarrow 2m$. le
tournoi $T-0$ \'etant ind\'ecomposable, il existe deux sommets $k$
et $l$ de $T-0$ tels que $k \longrightarrow 1$ et $2m
\longrightarrow l$. D'apr\`es ce qui pr\'ec\`ede, $k$ est impair
et $l$ est pair. Comme $k \longrightarrow 1$ et $k$ est impair
(resp. $2m \longrightarrow l$ et $l$ est pair), alors $\{i \in S
-\{1\} : i$ est impair\} $\longrightarrow 1$ (resp.  $2m
\longrightarrow \{i \in S-\{0, 2m\} : i$ est pair\}). Ainsi, $\{i
\in S - \{1\} : i$ est impair\} $\longrightarrow \{1, 2m\}
\longrightarrow \{i \in S-\{0, 2m\} : i$ est pair\}. En
particulier, $\{1, 2m\}$ est un intervalle non trivial du tournoi
ind\'ecomposable $T-0$. Contradiction.
\end {itemize}

{\hspace*{\fill}$\Box$\medskip}

\section{Preuve du th\'eor\`eme 1.3}

Nous montrons d'abord que les tournois introduits dans le
th\'eor\`eme $3.1$ sont $(-1)$-critiques et deux \`a deux non
isomorphes.

\begin{proposition} Soient $n$ et $k$ deux entiers tels que $n \geq
3$ et $k \in \{1, \cdots, n-2\}$. Les tournois $E_{2n+1}^{2k+1}$,
$F_{2n+1}^{2k+1}$, $G_{2n+1}^{2k+1}$ et $H_{2n+1}^{2k+1}$ sont des
tournois $(-1)$-critiques dont l'unique sommet non critique est
$2k+1$.
\end{proposition}

\noindent{\em Preuve\/. }

V\'erifions d'abord que pour $W = E_{2n+1}^{2k+1}$,
$F_{2n+1}^{2k+1}$, $G_{2n+1}^{2k+1}$ ou  $H_{2n+1}^{2k+1}$, on a:
\begin{equation}
\forall  i \in S(W) - \{2k+1\},\mbox{ }  W - \{i\} \mbox{ est
d\'ecomposable}.
\end{equation}

En effet, pour $i =0$ (resp. $i = 2n$), $\{2, \cdots, 2n\}$ (resp.
$\{0, \cdots, 2n-2\}$) est un intervalle non trivial de $W -i$.
Pour $i \in \{1, \cdots, 2k-2\} \cup \{2k\} \cup \{2k+2, \cdots,
2n-2\}$, $\{i-1, i+1\}$ est un intervalle non trivial de $W - i$.
Pour $i = 2k - 1$, $\{i-1, i+1\}$ est un intervalle non trivial de
chacun des tournois $E_{2n+1}^{2k+1} - i$, $F_{2n+1}^{2k+1} - i$
et $G_{2n+1}^{2k+1} - i$, et $\{i-1, i+2\}$ est un intervalle non
trivial de $H_{2n+1}^{2k+1} - i$. Pour $i = 2n-1$, $\{i-1, i+1\}$
(resp. $\{0, \cdots, 2n-3\} \cup \{2n\}$) est un intervalle non
trivial de chacun des tournois $E_{2n+1}^{2k+1} - i$ et
$F_{2n+1}^{2k+1} - i$ (resp. $G_{2n+1}^{2k+1} - i$ et
$H_{2n+1}^{2k+1} - i$).

Fixons maintenant un entier $k \geq 1$, et montrons la proposition
par r\'ecurrence sur l'entier $n \geq k+2$. Commen\c{c}ons par
examiner le cas o\`u $n = k+2$. Nos tournois sont maintenant
d\'efinis sur $\{0, \cdots, 2k+4\}$, on pose $X = \{0, \cdots,
2k+4\} - \{2k+1, 2k+2\}$. Les tournois $E_{2k+5}^{2k+1}(X)$ ,
$F_{2k+5}^{2k+1}(X)$, $G_{2k+5}^{2k+1}(X)$ et $H_{2k+5}^{2k+1}(X)$
sont ind\'ecomposables. En effet, d'une part $E_{2k+5}^{2k+1}(X)
\simeq H_{2k+5}^{2k+1}(X) \simeq V_{2k+3}$ (les $(2k+2)$-cha\^ines
respectives des tournois $E_{2k+5}^{2k+1}(X)$ et
$H_{2k+5}^{2k+1}(X)$ sont $0 < \cdots < 2k < 2k+3$ et $0 < \cdots
<2k-1 < 2k+3 < 2k+4$); d'autre part, $F_{2k+5}^{2k+1}(X) =
G_{2k+5}^{2k+1}(X) \simeq U_{2k+3}$ (un isomorphisme sur
$U_{2k+3}$ fixe chacun des sommets de $\{0, \cdots, 2k\}$, envoie
$2k+3$ sur $2k+1$ et $2k+4$ sur $2k+2$). Soit $T =
E_{2k+5}^{2k+1}$, $F_{2k+5}^{2k+1}$, $G_{2k+5}^{2k+1}$ ou
$H_{2k+5}^{2k+1}$. Dans le tournoi $T$, $2k+1 \in X(2k+3)$ et,
comme $1 \longrightarrow 2k+2 \longrightarrow 2k+3$, alors $2k+2
\notin [X]$. Montrons que $2k+2 \in Ext(X)$, ce qui signifie que :

\begin{equation}
T - \{2k+1\} \mbox{ est ind\'ecomposable }.
\end{equation}
Supposons par l'absurde qu'il existe $u \in X$ tel que $2k+2 \in
X(u)$. Pour $T = E_{2k+5}^{2k+1}$ ou $F_{2k+5}^{2k+1}$, $2k+2
\longrightarrow \{0, 2k+4\}$ et, comme il n'existe aucun sommet $x
\in X$ v\'erifiant $x \longrightarrow \{0, 2k +4\}$, alors $u \in
\{0, 2k+4\}$, ce qui contredit les faits que $0 \longrightarrow 1
\longrightarrow 2k+2$ et $2k+2 \longrightarrow 2k+3
\longrightarrow 2k+4$. Pour $T = G_{2k+5}^{2k+1}$ (resp.
$H_{2k+5}^{2k+1}$), $\{0, \cdots, 2k\} \longrightarrow 2k+2$ et
$T(\{0, \cdots, 2k\}) = U_{2k+1}$ (resp. $V_{2k+1}$), en
particulier $T(\{0, \cdots, 2k\})$ est ind\'ecomposable. Il
s'ensuit que $u \notin \{0, \cdots,2k\}$, autrement, $ \{0,
\cdots, 2k\} - \{u\} \longrightarrow u$, ce qui contredit
l'ind\'ecomposabilit\'e de $T(\{0, \cdots, 2k\})$. Ainsi, $u \in
\{2k+3, 2k+4\}$, contredisant le fait que $2k+2 \longrightarrow
2k+3 \longrightarrow 2k+4 \longrightarrow 2k+2$.

Comme $2k+1 \in X(2k+3)$ et $2k+2 \in Ext(X)$ et $2k+1
\longrightarrow 2k+2 \longrightarrow 2k+3$, alors, d'apr\`es le
lemme 2.2, le tournoi $T$ est ind\'ecomposable et, d'apr\`es
$(4.1)$ et $(4.2)$, $2k+1$ est l'unique sommet non critique de
$T$.

\`A pr\'esent, soit un entier $n > k + 2$. On pose $X_{n} = \{0,
\cdots, 2n\} - \{2k+2, 2k+3\}$ et $X'_{n} = X_{n} - \{2k +1\}$.
Soit $D = E$, $F$, $G$ ou $H$. L'application $f$ de $X_{n}$ sur
$\{0, \cdots, 2n-2\}$ d\'efinie par $f(i) = i$ (resp. $f(i) =
i-2$) si $i \in \{0, \cdots, 2k+1\}$ (resp. $i \in \{2k+4, \cdots,
2n\}$), est un isomorphisme de $D_{2n+1}^{2k+1}(X_{n})$ sur
$D_{2n-1}^{2k+1}$, qui fixe le sommet $2k+1$. Il s'ensuit, en
appliquant l'hypoth\`ese de r\'ecurrence, que
$D_{2n+1}^{2k+1}(X_{n})$ est un tournoi $(-1)$-critique dont
l'unique sommet non critique est $2k+1$. Le tournoi
$D_{2n+1}^{2k+1}(X'_{n})$ est ainsi ind\'ecomposable. On a bien
$2k+2 \in X_{n}(2k+4)$ et $2k+3 \in X_{n}(2k+1)$, en particulier,
$2k+2 \in X'_{n}(2k+4)$ et $D_{2n+1}^{2k+1}(X'_{n} \cup \{2k+3\})
\simeq D_{2n+1}^{2k+1}(X_{n}) \simeq D_{2n-1}^{2k+1}$. Le tournoi
$D_{2n+1}^{2k+1}(X'_{n} \cup \{2k+3\})$ \'etant alors
ind\'ecomposable, $2k+3 \in Ext(X'_{n})$. Comme de plus, $2k+2
\longrightarrow 2k+3 \longrightarrow 2k+4$, alors, d'apr\`es le
lemme 2.2, les tournois $D_{2n+1}^{2k+1}$ et $D_{2n+1}^{2k+1} -
\{2k+1\}$ sont ind\'ecomposables. Avec $(4.1)$, $2k+1$ est
l'unique sommet non critique de $D_{2n+1}^{2k+1}$.

{\hspace*{\fill}$\Box$\medskip}

\begin{corollaire}
Pour tout $n \geq 3$, les $6(n-2)$ tournois de la classe
$\mathcal{D}_{2n+1} = \mathcal{E}_{2n+1} \cup \mathcal{F}_{2n+1}
\cup \mathcal{F}_{2n+1}^{\star} \cup \mathcal{G}_{2n+1} \cup
\mathcal{G}_{2n+1}^{\star} \cup \mathcal{H}_{2n+1}$ sont deux \`a
deux non isomorphes.

\end{corollaire}
\noindent{\em Preuve\/. }

Soient $T$ et $T'$ deux tournois isomorphes de la classe
$\mathcal{D}_{2n+1}$. D'apr\`es la proposition 4.1, $T$ et $T'$
sont (-1)-critiques. D\'esignons alors par $a$ et $a'$ leurs
sommets non critiques respectifs. Par construction des
diff\'erentes classes, le tournoi $T$ est dans la classe
$\mathcal{E}_{2n+1}$ (resp. $\mathcal{F}_{2n+1}$,
$\mathcal{F}_{2n+1}^{\star}$, $\mathcal{G}_{2n+1}$,
$\mathcal{G}_{2n+1}^{\star}$ $\mathcal{H}_{2n+1}$) si et seulement
si $T(V_{T}^{+}(a))$ et $T(V_{T}^{-}(a))$ sont des cha\^ines
(resp. $T(V_{T}^{+}(a))$ est une cha\^ine et $T(V_{T}^{-}(a))$
n'est pas une cha\^ine, $T(V_{T}^{+}(a))$ n'est pas une cha\^ine
et $T(V_{T}^{-}(a))$ est une cha\^ine, il existe un unique sommet
de $V_{T}^{+}(a)$ qui domine au moins deux sommets de
$V_{T}^{-}(a)$, il existe un unique sommet de $V_{T}^{-}(a)$ qui
soit domin\'e par au moins deux sommet de $V_{T}^{+}(a)$, $\mid
A(T) \cap (V_{T}^{+}(a) \times V_{T}^{-}(a)) \mid = 1$). Il
s'ensuit que si $T$ est un tournoi de la classe
$\mathcal{E}_{2n+1}$ (resp. $\mathcal{F}_{2n+1}$,
$\mathcal{F}_{2n+1}^{\star}$, $\mathcal{G}_{2n+1}$
$\mathcal{G}_{2n+1}^{\star}$, $\mathcal{H}_{2n+1}$), il en est de
m\^eme pour le tournoi $T'$. Si $T$ et $T'$ sont dans la m\^eme
classe $\mathcal{E}_{2n+1}$ (resp. $\mathcal{F}_{2n+1}$,
$\mathcal{G}_{2n+1}$, $\mathcal{H}_{2n+1}$), alors il existe deux
entiers $k$, $k' \in \{0, \cdots, n-2\}$ tels que $T =
E_{2n+1}^{2k+1}$ (resp. $F_{2n+1}^{2k+1}$, $G_{2n+1}^{2k+1}$,
$H_{2n+1}^{2k+1}$) et $T' = E_{2n+1}^{2k'+1}$ (resp.
$F_{2n+1}^{2k'+1}$, $G_{2n+1}^{2k'+1}$, $H_{2n+1}^{2k'+1}$).
D'apr\`es la proposition 4.1, $a = 2k+1$ et $a' = 2k' + 1$. Ainsi,
$\mid V_{T}^{-}(a) \mid = 2k +1$ et $\mid V_{T'}^{-}(a') \mid =
2k' +1$. Or, comme un isomorphisme de $T$ sur $T'$ envoie $a$ sur
$a'$, on a $\mid V_{T}^{-}(a) \mid = \mid V_{T'}^{-}(a') \mid$. Il
s'ensuit que $k = k'$ et donc $T = T'$. Si enfin $T$ et $T'$ sont
dans la m\^eme classe $\mathcal{F}_{2n+1}^{\star}$ (resp.
$\mathcal{G}_{2n+1}^{\star}$), alors  $T^{\star}$ et $T'^{\star}$
sont dans la m\^eme classe $\mathcal{F}_{2n+1}$ (resp.
$\mathcal{G}_{2n+1}$) de sorte que, d'apr\`es ce qui pr\'ec\`ede,
$T^{\star} = T'^{\star}$ et donc $T = T'$.

{\hspace*{\fill}$\Box$\medskip}

Il est commode d'introduire les deux lemmes suivants, avant
d'entamer la preuve du th\'eor\`eme.

\begin{lemme}

Soit $T = (S, A)$ un tournoi (-1)-critique, avec
$I(T)(\mathcal{C}) = P_{m+1}$ o\`u $\mathcal{C} = \{0, \cdots,
m\}$ est la composante connexe, non r\'eduite \`a un singleton, de
$I(T)$. On d\'esigne par $a$ le sommet non critique de $T$ et on
pose $A_{I}^{+} = \{i \in \mathcal{C}: i>a$ et $i$ est impair\},
$A_{I}^{-} = \{i \in \mathcal{C}: i<a$ et $i$ est impair\},
$A_{P}^{+} = \{i \in \mathcal{C}: i>a$ et $i$ est pair\} et
$A_{P}^{-} = \{i \in \mathcal{C}: i<a$ et $i$ est pair\}. Si $W =
T(A_{P}^{+})$, $T(A_{P}^{-})$, $T(A_{I}^{+})$ ou $T(A_{I}^{-})$,
alors $W$ ou $W^{\star}$ est l'ordre total usuel sur l'ensemble
$S(W)$ des sommets de $W$.

\end{lemme}

\noindent{\em Preuve\/. }

Il suffit de consid\'erer le cas o\`u $\mid S(W) \mid \geq 3$. On
pose alors $S(W) =\{w_{1}, \cdots, w_{q}\}$, o\`u $q \geq 3$ et
$w_{1} < \cdots < w_{q}$. Soit $k \in \{2, \cdots, q\}$. On a
$w_{k-1} +1 \in \mathcal{C} - \{0, m\}$ et donc $\mid
V_{I(T)}(w_{k-1} + ~ 1)\mid = 2$. En utilisant le lemme 3.1,
$\{w_{k-1}, w_{k}\}$ est un intervalle de $T- \{w_{k-1}+1\}$ et,
comme $w_{k-1} +1 \notin S(W)$, alors $\{w_{k-1}, w_{k}\}$ est un
intervalle de $W$. Il s'ensuit que pour tous $i < j \in
\{1,\cdots, q\}$, $(w_{i},w_{j}) \equiv \cdots \equiv
(w_{1},w_{j}) \equiv (w_{1},w_{j-1}) \equiv \cdots \equiv
(w_{1},w_{2})$. Il en d\'ecoule que si $w_{1} \longrightarrow
w_{2}$ (resp. $w_{2} \longrightarrow w_{1}$ ), alors $W$ (resp.
$W^{\star}$) est l'ordre total usuel sur $S(W)$.

{\hspace*{\fill}$\Box$\medskip}

\begin{lemme}
Soit $T = (S, A)$ un tournoi (-1)-critique, avec $S = \{0, \cdots,
2n\}$ et $I(T)(\mathcal{C}) = P_{m+1}$ o\`u $\mathcal{C} = \{0,
\cdots, m\}$ est la composante connexe, non r\'eduite \`a un
singleton, de $I(T)$. Alors $m \geq 3$, et si $0 \longrightarrow
1$, alors $V_{T}^{+}(1) =  \{2, \cdots, 2n\}$ et $V_{T}^{+}(m-1) =
\{m\}$. De plus, en d\'esignant par $a$ le sommet non critique de
$T$, les assertions suivantes sont v\'erifi\'ees.
\begin{enumerate}
\item Si $a$ est impair, alors pour tout sommet impair $i$, on a:
\begin{itemize}

\item Si $i \in \{1, \cdots, m\}$, alors $V_{T}^{+}(i) = \{i+1,
\cdots, 2n\}$.
\item Si $i \in \{1, \cdots, m - a\}$, alors $V_{T}^{+}(m-i) = \{m -i + 1,
\cdots, m\}$.
\end{itemize}
En particulier, ou bien $S - \mathcal{C} = \emptyset$, ou bien $m$
est impair.
\item Si $a$ et $m$ sont pairs, alors pour tout sommet impair $i$,
on a:
\begin{itemize}
\item Si $i \in \{1, \cdots, a - 1\}$, alors $V_{T}^{+}(i) = \{i+1,
\cdots, 2n\}$.
\item Si $i \in \{1, \cdots, m-a-1\}$, alors $V_{T}^{+}(m-i) = \{m-i+1,
\cdots, m\}$.

\end{itemize}
En particulier, $S - \mathcal{C} \neq \emptyset$.
\end{enumerate}
\end{lemme}

\noindent{\em Preuve\/. }

Notons d'abord que, d'apr\`es les lemmes 3.5 et 3.4,
$I(T)(\mathcal{C})$ est un chemin et $a \in \{1, \cdots, m-1\}$.
On a bien $m \geq 3$, autrement, $m = 2$ et $a = 1$ et dans ce
cas, d'apr\`es le lemme 3.1, $S - \{0, 1\}$ et $S - \{1, 2\}$ sont
des intervalles respectifs des tournois $T-0$ et $T - 2$, de sorte
que $1 \sim S - \{0, 1\}$ et $1 \sim S - \{1, 2\}$. Comme de plus,
$(S - \{0, 1\}) \cap (S - \{1, 2\}) \neq \emptyset$, alors $1 \sim
S -\{1\}$, ce qui contredit l'ind\'ecomposabilit\'e de $T$.
D'apr\`es le lemme 3.1, $1 \sim S - \{0, 1\}$ et $m-1 \sim S -
\{m-1, m\}$ et, comme $T$ est ind\'ecomposable avec $0
\longrightarrow 1$, alors $V_{T}^{+}(1) = \{2, \cdots, 2n\}$, en
particulier, $1 \longrightarrow m-1$ et donc $V_{T}^{+}(m - 1) =
\{m\}$. Supposons, \`a pr\'esent, que $a$ est impair et montrons,
par r\'ecurrence finie, l'assertion $1$. D'apr\`es ce qui
pr\'ec\`ede, l'assertion est v\'erifi\'ee pour $i = 1$. Si $i$ est
un sommet impair avec $3 \leq i \leq m$ (resp. $3 \leq i \leq m -
a $), alors, par hypoth\`ese de r\'ecurrence, $V_{T}^{+}(i-2) =
\{i-1, \cdots, 2n\}$ (resp. $V_{T}^{+}(m-i+2) = \{m-i+3, \cdots,
m\}$). Dans ce cas, $i -1$ (resp. $m-i+1$) est un sommet critique
de $T$ avec $V_{I(T)}(i-1) = \{i, i-2\}$ (resp. $V_{I(T)}(m-i+1) =
\{m-i, m-i+2\}$) de sorte que, d'apr\`es le lemme 3.1, $\{i,
i-2\}$ (resp. $\{m-i, m-i+2\}$) est un intervalle de $T-\{i-1\}$
(resp. $T - \{m-i+1\}$), sans \^etre un intervalle de $T$. Il
s'ensuit, en utilisant l'hypoth\`ese de r\'ecurrence, que
$V_{T}^{+}(i) = \{i+1, \cdots, 2n\}$ (resp. $V_{T}^{+}(m-i) =
\{m-i+1, \cdots, m\}$). En particulier, si $m$ est pair, alors
$V_{T}^{+}(m-1) = \{m\} = \{m, \cdots, 2n\}$, de sorte que $m =
2n$ et donc $S - \mathcal{C} = \emptyset$.

Supposons maintenant que $a$ et $m$ sont pairs et montrons, par
r\'ecurrence finie, l'assertion $2$. Celle-ci est v\'erifi\'ee
pour $i = 1$. Si $i$ est un sommet impair avec $3 \leq i \leq a-1$
(resp. $3 \leq i \leq m-a-1$), alors, par hypoth\`ese de
r\'ecurrence, $V_{T}^{+}(i-2) = \{i-1, \cdots, 2n\}$ (resp.
$V_{T}^{+}(m-~i+2) = \{m-i+3, \cdots, m\}$). Dans ce cas, $i-1$
(resp. $m-i+1$) est un sommet critique de $T$ et, en utilisant le
lemme 3.1, $\{i, i-2\}$ (resp. $\{m-i, m-i+2\}$) est un intervalle
de $T - \{i-1\}$ (resp. $T - \{m-i+1\}$), sans \^etre un
intervalle de $T$. Il s'ensuit, en utilisant l'hypoth\`ese de
r\'ecurrence, que $V_{T}^{+}(i) = \{i+1, \cdots, 2n\}$ (resp.
$V_{T}^{+}(m-i) = \{m-i+1, \cdots, m\}$). En particulier, si $S -
\mathcal{C} = \emptyset$, c'est \`a dire $m = 2n$, alors
$V_{T-a}^{+}(a-1) = \{a+1, \cdots, 2n\}$ et $V_{T-a}^{+}(a+1) =
\{a+2, \cdots, 2n\}$, de sorte que $\{a-1, a+1\}$ est un
intervalle non trivial du tournoi ind\'ecomposable $T-a$.
Contradiction.

{\hspace*{\fill}$\Box$\medskip}

Maintenant, nous pr\'esentons une preuve du th\'eor\`eme $1.3$,
qui repose sur une construction des tournois $(-1)$-critiques \`a
partir des diff\'erents graphes d'ind\'ecomposabilit\'es possibles
pour de tels tournois.

{\bfseries Th\'eor\`eme 1.3} \begin{em}\`A un isomorphisme pr\`es,
les tournois $(-1)$-critiques sont les tournois $E_{2n+1}^{2k+1}$,
$F_{2n+1}^{2k+1}$, $(F_{2n+1}^{2k+1})^{\star}$, $G_{2n+1}^{2k+1}$,
$(G_{2n+1}^{2k+1})^{\star}$ et $H_{2n+1}^{2k+1}$, o\`u $n \geq 3$
et $1 \leq k \leq n-2$. De plus, le sommet $2k+1$ est l'unique
sommet non critique de chacun de ces tournois.\end{em}

\noindent{\em Preuve\/. }

Soit $T = (S, A)$ un tournoi $(-1)$-critique. D'apr\`es la
proposition 4.1 et le corollaire 4.2, il suffit de montrer que $T$
est dans la classe $\mathcal{E}_{2n+1} \cup \mathcal{F}_{2n+1}
\cup \mathcal{F}_{2n+1}^{\star} \cup \mathcal{G}_{2n+1} \cup
\mathcal{G}_{2n+1}^{\star} \cup \mathcal{H}_{2n+1}$, pour un
entier $n \geq 3$. On d\'esigne par $a$ l'unique sommet non
critique de $T$ et par $\mathcal{C}$ la composante connexe de
$I(T)$, non r\'eduite \`a un singleton. On pose $I(T)(\mathcal{C})
= P_{m +1}$ et $S = \{0, \cdots, 2n\}$, o\`u $ n \geq 3$ et $m
\leq 2n$. D'apr\`es le lemme 3.4, $a \in \{1, \cdots, m-1\}$. On
pose $A_{I}^{+} = \{i \in \mathcal{C}: i>a$ et $i$ est impair\},
$A_{I}^{-} = \{i \in \mathcal{C}: i<a$ et $i$ est impair\},
$A_{P}^{+} = \{i \in \mathcal{C}: i>a$ et $i$ est pair\},
$A_{P}^{-} = \{i \in \mathcal{C}: i<a$ et $i$ est pair\}, $A^{+} =
A_{P}^{+} \cup A_{I}^{+}$ et $A^{-} = A_{P}^{-} \cup A_{I}^{-}$.
Quitte \`a remplacer $T$ par $T^{\star}$, on peut supposer que $0
\longrightarrow 1$.

 Supposons d'abord que $a$ est
impair. Soit $(i, j) \in A_{P}^{-} \times A_{P}^{+}$. D'une part,
si $i \neq a-1$ alors, d'apr\`es le lemme 3.1, $\{i, i+2\}$ est un
intervalle de $T - \{i+1\}$, de sorte que $(i, j) \equiv (i+2,
j)$. D'autre part, si $j \leq m-2$, alors $\{j, j+2\}$ est un
intervalle de $T - \{j+1\}$, de sorte que $(i, j) \equiv (i,
j+2)$. Il s'ensuit que :
\begin{equation}
A_{P}^{-} \sim A_{P}^{+}.
\end{equation}

 Le lemme 4.4 nous am\`ene \`a distinguer les deux cas suivants.
\begin{itemize}
\item $S - \mathcal{C} = \emptyset$. Dans ce cas $m =2n$, et le
lemme 4.4 donne que :

\begin{equation}
\mbox{ pour tout sommet impair } i \mbox{ de } T, \mbox{ }
V_{T}^{+}(i) = \{i+1, \cdots, 2n\}.
\end{equation}

N\'ecessairement $a+1 \longrightarrow 0$, sinon, d'apr\`es
$(4.3)$, $A_{P}^{-} \longrightarrow A_{P}^{+}$ et, en utilisant
$(4.3)$, on obtient que $A^{-} \longrightarrow A^{+}$, ce qui
contredit l'ind\'ecomposabilit\'e de $T-a$. Il s'ensuit, en
utilisant (3), que :
\begin{equation}
A_{P}^{+} \longrightarrow A_{P}^{-}.
\end{equation}

On distingue les cas suivants.
\begin{itemize}

\item $0 \longrightarrow 2$. Comme $a+1 \longrightarrow 0 \longrightarrow
2$, et $a$ est impair, alors $a \geq 3$. Le lemme 4.3 donne que
$T(A_{P}^{-})$ est l'ordre total usuel sur $A_{P}^{-}$. En
utilisant de plus la relation $(4.4)$, on obtient:

\begin{equation}
T(A^{-}) = L_{a}.
\end{equation}

 Notons que $a \leq 2n-3$, autrement, $a = 2n-1$ et, en utilisant $(4.4)$, $(4.5)$ et $(4.6)$, $T =
V_{2n+1}$, en particulier, $T$ est un tournoi critique,
contradiction. Comme de plus  $a \geq 3$, alors $3 \leq a \leq
2n-3$. On pose alors $a = 2k + 1$, avec $k \in \{1, \cdots,
n-~2\}$.

Si $a+1 \longrightarrow a+3$, le lemme 4.3 donne que
$T(A_{P}^{+})$ est l'ordre total usuel sur $A_{P}^{+}$ et, en
utilisant de plus $(4.4)$, $(4.5)$ et $(4.6)$, on obtient que $T =
E_{2n+1}^{2k+1}$.

Si au contraire, $a+3 \longrightarrow a+1$, alors $T(A_{P}^{+})$
est, encore par le lemme 4.3, le tournoi dual de l'ordre total
usuel sur $A_{P}^{+}$ et, en utilisant encore $(4.4)$, $(4.5)$ et
$(4.6)$, on trouve que $T \simeq (F_{2n+1}^{2k'+1})^{\star}$, o\`u
$k' = n-k-1 \in \{1, \cdots, n-2\}$ (un isomorphisme $f$ de $T$
sur $(F_{2n+1}^{2k'+1})^{\star}$ est d\'efini par : pour tout $i
\in \{0, \cdots, 2n\}$, $f(i) = 2n-i$).

\item $2 \longrightarrow 0$.
Dans ce cas, d'apr\`es le lemme 4.3, $T(A_{P}^{-})$ est le tournoi
dual de l'ordre total usuel sur $A_{P}^{-}$. Notons que $a \leq
2n-3$, autrement, $a = 2n-1$ et, en utilisant de plus $(4.4)$ et
$(4.5)$, on obtient que $T = U_{2n+1}$, en particulier, $T$ est un
tournoi critique, contradiction.

N\'ecessairement $a+1 \longrightarrow a+3$, sinon, par le lemme
4.3, $T(A_{P}^{+})$ est le tournoi dual de l'ordre total usuel sur
$A_{P}^{+}$ et, en utilisant, en outre, $(4.4)$ et $(4.5)$, on
obtient que $T = U_{2n+1}$, en particulier $T$ est critique,
contradiction. Il s'ensuit, encore par le lemme 4.3, que
$T(A_{P}^{+})$ est l'ordre total usuel sur $A_{P}^{+}$. En
particulier, $a \geq 3$, autrement, $a = 1$ et en utilisant, de
plus, $(4.4)$ et $(4.5)$, on obtient que $T \simeq V_{2n+1}$ (la
$2n$-cha\^ine de $T$ \'etant: $1 < 2 < \cdots <2n$). Comme de plus
$a \leq 2n-3$, alors $3 \leq a \leq 2n-3$. On pose alors $a = 2k +
1$, avec $k \in \{1, \cdots, n-2\}$. En utilisant $(4.4)$, $(4.5)$
et le fait que $T(A_{P}^{+})$ (resp. $T(A_{P}^{-})$) est l'ordre
total usuel sur $T(A_{P}^{+})$ (resp. le dual de l'ordre total
usuel sur $T(A_{P}^{-})$), on obtient que $T = F_{2n+1}^{2k+1}$.

\end{itemize}
\item m est impair. Dans ce cas $S - \mathcal{C} \neq \emptyset$,
et le lemme 4.4 nous donne les deux faits suivants.
\begin{equation}
\mbox{ pour tout sommet impair } i \in \{1, \cdots, m\},\
V_{T}^{+}(i) = \{i+1, \cdots, 2n\}.
\end {equation}

\begin{equation}
T(A_{P}^{+}) \mbox{ est l'ordre total usuel sur } A_{P}^{+},
\mbox{ avec } A^{-} \cup (S - \mathcal{C}) \longrightarrow
A_{P}^{+} .
\end{equation}

Soit $\mu \in S - \mathcal{C}$ et soit $i \in A_{P}^{-} - \{0\}$.
Comme $\{i, i-2\}$ est un intervalle de $T - \{i-1\}$ et $\mu \neq
i-1$, alors $(i, \mu) \equiv (i-2, \mu) \equiv \cdots \equiv (0,
\mu)$. Il s'ensuit que :

\begin{equation}
\mbox{ pour tout } \mu \in S - \mathcal{C},\mbox{  } \mu \sim
A_{P}^{-}.
\end{equation}

Il existe alors $\lambda \in S - \mathcal{C}$ tel que $\lambda
\longrightarrow 0$. Autrement, $0 \longrightarrow S - \mathcal{C}$
et, en utilisant $(4.9)$, on obtient que $A_{P}^{-}
\longrightarrow S- \mathcal{C}$. Il s'ensuit, en utilisant de plus
$(4.7)$ et $(4.8)$, que $A^{-} \longrightarrow S - A^{-}$, ce qui
contredit l'ind\'ecomposabilit\'e de $T$. Ainsi, de nouveau par
$(4.9)$, $\lambda \longrightarrow A_{P}^{-}$. Supposons par
l'absurde que $0 \longrightarrow 2$. En utilisant $(4.7)$, $(4.8)$
et le lemme 4.3, on obtient que $T(\mathcal{C}) = L_{m+1}$ et on
d\'eduit, \`a l'aide de $(4.7)$, $(4.8)$, $(4.9)$ et le fait que
$\lambda \longrightarrow A_{P}^{-}$, que $T(\mathcal{C} \cup
\{\lambda\}) \simeq V_{m+2}$. Comme le tournoi $T$ n'est pas
critique, $S - \mathcal{C} \neq \{\lambda\}$, et puisque $\mid S -
\mathcal{C}\mid$ est impair, alors, par une suite finie
d'applications du corollaire 2.3 au sous-tournoi ind\'ecomposable
$T(\mathcal{C} \cup \{\lambda\})$ de $T$, on obtient deux sommets
distincts $x$ et $y$ de $S - (\mathcal{C} \cup \{\lambda\})$ tels
que $T -\{x, y\}$ est ind\'ecomposable. Ceci contredit le fait que
$\{x, y\}$ n'est pas une ar\^ete de $I(T)$. Ainsi, forc\'ement $2
\longrightarrow 0$ et, d'apr\`es le lemme 4.3, $T(A_{P}^{-})$ est
le tournoi dual de l'ordre total usuel sur $A_{P}^{-}$. De plus,
$a \neq 1$, autrement, $A^{-} = \{0\}$ et $2 \in A_{P}^{+}$ de
sorte que, d'apr\`es $(4.8)$, $0 \longrightarrow 2$,
contradiction. On pose alors $n' = \frac{m+1}{2}$ et $a = 2k +1$,
o\`u $k \in \{1, \cdots, n'-2\}$. Le tournoi $T(\mathcal{C} \cup
\{\lambda\})$ est isomorphe \`a $G_{2n'+1}^{2k+1}$ (un
isomorphisme de $T(\mathcal{C} \cup \{\lambda\})$ sur
$G_{2n'+1}^{2k+1}$ envoie $\lambda$ sur $2n'$ et fixe chaque
sommet de $\mathcal{C}$). En particulier, par la proposition
$4.1$, le tournoi $T(\mathcal{C} \cup \{\lambda\})$ est
ind\'ecomposable. Il s'ensuit que $S- \mathcal{C} = \{\lambda\}$.
Autrement, $\mid S - \mathcal{C}\mid$ est impair avec $\mid S -
\mathcal{C} \mid \geq 3$ et, par une suite finie d'applications le
corollaire 2.3 au sous-tournoi ind\'ecomposable $T(\mathcal{C}
\cup \{\lambda\})$ de $T$, on obtient deux sommets distincts $x$
et $y$ de $S - \mathcal{C}$ tels que $T - \{x, y\}$ est
ind\'ecomposable. Ceci contredit le fait que $x$ est un sommet
isol\'e de $I(T)$. On conclut que $n' = n$ et que $T =
T(\mathcal{C} \cup \{\lambda\}) \simeq G_{2n+1}^{2k+1}$.
\end{itemize}
Il reste \`a examinier le cas o\`u $a$ est pair. Soit, dans ce
cas, $\sigma$ la permutation de $S$ fixant chaque sommet de $S -
\mathcal{C}$ et telle que pour tout $i \in \mathcal{C}$,
$\sigma(i) = m-i$. L'application $\sigma$ est un isomorphisme de
$T$ sur un tournoi $T'$. Si m est impair, alors le sommet non
critique de $T'$ est le sommet impair $\sigma(a) = m-a$. On se
ram\`ene ainsi au cas pr\'ecedent et on d\'eduit que $T'$, et par
suite $T$, est un tournoi de la classe $\mathcal{E}_{2n+1} \cup
\mathcal{F}_{2n+1} \cup (\mathcal{F}_{2n+1})^{\star} \cup
\mathcal{G}_{2n+1} \cup (\mathcal{G}_{2n+1})^{\star}$, pour un
entier $n \geq 3$. Supposons alors que $m$ est pair. D'apr\`es le
lemme 4.4, d'une part $m \geq 4$, d'autre part, le tournoi $T$ est
d\'etermin\'e par $T( S - (A_{I}^{-} \cup A_{I}^{+}))$. Posons $R
= T(\mathcal{C}_{P})$, o\`u $\mathcal{C}_{P}$ est l'ensemble des
sommets pairs de $\mathcal{C}$. On a bien $R$ ou $R^{\star}$ est
l'ordre total usuel sur $\mathcal{C}_{P}$. En effet, si $i$ et $j$
sont deux sommets de $\mathcal{C}_{P}$, avec $0 < i < j$, alors,
en utilisant le lemme 3.1, $(i, j) \equiv (i-2, j) \equiv \cdots
(0, j) \equiv (0, j-2)\equiv \cdots \equiv(0, 2)$. Il s'ensuit que
si $0 \longrightarrow 2$ (resp. $2 \longrightarrow 0$), alors $R$
(resp. $R^{\star}$) est l'ordre total usuel sur $C_{P}$. Supposons
par l'absurde que $2 \longrightarrow 0$. D'apr\`es ce qui
pr\'ec\`ede, $R$ est le tournoi dual de l'ordre total usuel sur
$\mathcal{C}_{P}$ et, en utilisant de plus le lemme 4.4,
$T(\mathcal{C}) = U_{m+1}$. Comme de plus $\mid S - \mathcal{C}
\mid$ est pair avec $\mid S - \mathcal{C} \mid \geq 2$, alors, une
suite finie d'applications du corollaire 2.3 au
sous-tournoi-ind\'ecomposable $T(\mathcal{C})$ de $T$, nous donne
deux sommets distincts $x$ et $y$ de $S - \mathcal{C}$ tels que $T
- \{x, y\}$ est ind\'ecomposable. Ceci contredit le fait que $x$
est un sommet isol\'e de $I(T)$. Ainsi, $0 \longrightarrow 2$, de
sorte que $R$ est l'ordre total usuel sur $\mathcal{C}_{P}$ et,
avec le lemme 4.4, $T(\mathcal{C})$ est l'ordre total usuel sur
$\mathcal{C}$.

Soit $\nu \in S - \mathcal{C}$ et soit $l \in \mathcal{C}_{P} -
\{0\}$. D'apr\`es le lemme 3.1, $\{l, l-2\}$ est un intervalle de
$T - \{l-1\}$ et, comme $\nu \neq l-1$, alors $(l, \nu) \equiv
(l-2, \nu) \equiv \cdots \equiv (0, \nu)$. Il s'ensuit que:

\begin{equation}
\mbox{ pour tout } \nu \in S - \mathcal{C},\mbox{  } \nu \sim
\mathcal{C}_{P}.
\end{equation}

On en d\'eduit que $V_{T}^{+}(0) \cap (S - \mathcal{C}) \neq
\emptyset$. Autrement, $(S - \mathcal{C}) \longrightarrow 0$ et,
d'apr\`es $(4.10)$, $S - \mathcal{C} \longrightarrow
\mathcal{C}_{P}$, de sorte qu'en utilisant de plus le lemme 4.4,
on obtient que $S - \mathcal{C}$ est un intervalle non trivial du
tournoi ind\'ecomposable $T$. Contradiction. On a aussi
$V_{T}^{-}(0) \cap (S - \mathcal{C}) \neq \emptyset$. Sinon, comme
$T(\mathcal{C})$ est l'ordre total usuel sur $\mathcal{C}$,
$V_{T}^{-}(0) = \emptyset$, ce qui contredit
l'ind\'ecomposabilit\'e de $T$. On pose $\Gamma^{+} = V_{T}^{+}(0)
\cap (S - \mathcal{C})$ et $\Gamma^{-} = V_{T}^{-}(0) \cap
(S-\mathcal{C})$. D'apr\`es $(4.10)$, pour tout $(u, v) \in
\Gamma^{-} \times \Gamma^{+}$, $u \longrightarrow \mathcal{C}_{P}
\longrightarrow v$. Il existe alors $(\alpha, \beta) \in
\Gamma^{-} \times \Gamma^{+}$ tel que $\beta \longrightarrow
\alpha$, autrement, $\Gamma^{-}  \cup A^{-} \cup \{a\}
\longrightarrow \Gamma^{+} \cup A^{+}$, ce qui contredit
l'ind\'ecomposabilit\'e de $T$. Il s'ensuit que $T(\mathcal{C}
\cup \{\alpha, \beta\}) \simeq H_{2n'+1}^{2k+1}$, o\`u $n' =
\frac{m+2}{2} \geq 3$ et $k = \frac{a}{2} \in \{1, \cdots, n'-2\}$
(un isomorphisme $\phi$ de $T(\mathcal{C} \cup \{\alpha, \beta\})$
sur $H_{2n'+1}^{2k+1}$ est d\'efini par $\phi(i) = i$ (resp.
$i+1$), si $i \in \{0, \cdots, 2k-1\}$ (resp. $i \in \{2k, \cdots,
m\}$), $\phi(\alpha) = 2k$ et $\phi(\beta) = 2n'$). Il s'ensuit
que $S- \mathcal{C} = \{\alpha, \beta\}$, et donc $n' = n$ et $T =
T(\mathcal{C} \cup \{\alpha, \beta\}) \simeq H_{2n+1}^{2k+1}$.
Autrement, $\mid S - \mathcal{C} \mid$ \'etant pair, par une suite
finie d'applications du corollaire 2.3 au sous-tournoi
ind\'ecomposable $T(\mathcal{C} \cup \{\alpha, \beta\})$ de $T$,
on obtient deux sommets distincts $x$ et $y$ de $S - \mathcal{C}$
tels que $T - \{x, y\}$ est ind\'ecomposable. Ceci contredit le
fait que $x$ est un sommet isol\'e de $I(T)$.
{\hspace*{\fill}$\Box$\medskip}

Cette preuve du th\'eor\`eme $1.3$, donne les tournois
$(-1)$-critiques avec leurs graphes d'ind\'ecomposabilit\'e. Nous
en d\'egageons la remarque suivante.

\begin{remarque}
Le graphe d'ind\'ecomposabilit\'e d'un tournoi (-1)-critique admet
au plus deux sommets isol\'es. Plus pr\'ecis\'ement, pour $n \geq
3$ et pour $1\leq k \leq n-2$, on a :
\begin{itemize}
\item $I(E_{2n+1}^{2k+1}) = I(F_{2n+1}^{2k+1}) = P_{2n+1}$.
\item $I(G_{2n+1}^{2k+1})$ est obtenu \`a partir de $P_{2n+1}$ en supprimant
l'ar\^ete $\{2n-1, 2n\}$.
\item $I(H_{2n+1}^{2k+1})$ est obtenu \`a partir de $P_{2n+1}$ en
supprimant les trois ar\^etes $\{2n-1, 2n\}$, $\{2k - 1, 2k\}$,
$\{2k, 2k+1\}$, et en ajoutant la paire $\{2k-1, 2k+1\}$ comme
nouvelle ar\^ete.
\end{itemize}
\end{remarque}

Notons que, d'apr\`es le lemme 3.4, si $a$ est le sommet non
critique d'un tournoi $(-1)$- critique $T$ , alors $T-a$ est un
tournoi $(-2)$-critique. Cela nous am\`ene \`a poser le probl\`eme
suivant.

\begin{probleme}
Caract\'eriser les tournois $(-2)$-critiques.
\end{probleme}

\end{document}